\documentclass[11 pt]{amsart}

\usepackage{amsmath}
\usepackage{amssymb}
\usepackage{amsfonts}
\usepackage{amsthm}
\usepackage{enumerate}

\begin{document}
\newcommand{\reminder}[1]{\centerline{\fbox{\parbox{12.5cm}{#1}}}}

\theoremstyle{plain}
\newtheorem{thm}{Theorem}[section]
\newtheorem{lemma}[thm]{Lemma}
\newtheorem{corollary}[thm]{Corollary}
\numberwithin{equation}{section}

\newtheorem{maintheorem}{Main Theorem}

\theoremstyle{definition}
\newtheorem{definition}{Definition}[section]
\newtheorem{remark}{Remark}[section]
\newtheorem{problem}{Problem}[section]

\numberwithin{equation}{section}

\def\intslash{\rlap{\kern  .32em $\mspace {.5mu}\backslash$ }\int}
\def\qsl{{\rlap{\kern  .32em $\mspace {.5mu}\backslash$ }\int_{Q_x}}}
\def\Re{\operatorname{Re\,}}
\def\Im{\operatorname{Im\,}}
\def\mx{{\max}}
\def\mn{{\min}}
\def\vth{\vartheta}

\def\rn{\rr^{n}}
\def\rr{\mathbb R}
\def\R{\mathbb R}
\def\Q{\mathcal Q}
\def\N{\mathbb N}
\def\complex{{\mathbb C}}
\def\norm#1{{ \left|  #1 \right| }}
\def\Norm#1{{ \left\|  #1 \right\| }}
\def\set#1{{ \left\{ #1 \right\} }}
\def\floor#1{{\lfloor #1 \rfloor }}
\def\emph#1{{\it #1 }}
\def\diam{{\text{\rm diam}}}
\def\osc{{\text{\rm osc}}}
\def\ffB{\mathcal B}
\def\itemize#1{\item"{#1}"}
\def\seq{\subseteq}
\def\Id{\text{\sl Id}}

\def\Ga{\Gamma}
\def\ga{\gamma}
\def\Th{\Theta}

\def\prd{{\text{\it prod}}}
\def\parab{{\text{\it parabolic}}}

\def\eg{{\it e.g. }}
\def\cf{{\it cf}}
\def\Rn{{\mathbb R^n}}
\def\Rd{{\mathbb R^d}}
\def\sgn{{\text{\rm sign }}}
\def\rank{{\text{\rm rank }}}
\def\corank{{\text{\rm corank }}}
\def\coker{{\text{\rm Coker }}}
\def\loc{{\text{\rm loc}}}
\def\spec{{\text{\rm spec}}}

\def\comp{{\text{\rm comp}}}

\def\Coi{{C^\infty_0}}
\def\dist{{\text{\rm dist}}}
\def\diag{{\text{\rm diag}}}
\def\supp{{\text{\rm supp }}}
\def\rad{{\text{\rm rad}}}
\def\Lip{{\text{\rm Lip}}}
\def\inn#1#2{\langle#1,#2\rangle}
\def\biginn#1#2{\big\langle#1,#2\big\rangle}
\def\rta{\rightarrow}
\def\lta{\leftarrow}
\def\noi{\noindent}
\def\lcontr{\rfloor}
\def\lco#1#2{{#1}\lcontr{#2}}
\def\lcoi#1#2{\imath({#1}){#2}}
\def\rco#1#2{{#1}\rcontr{#2}}
\def\bin#1#2{{\pmatrix {#1}\\{#2}\endpmatrix}}
\def\meas{{\text{\rm meas}}}

\def\card{\text{\rm card}}
\def\lc{\lesssim}
\def\gc{\gtrsim}
\def\pv{\text{\rm p.v.}}

%Greek letters
\def\alp{\alpha}             \def\Alp{\Alpha}
\def\bet{\beta}
\def\gam{\gamma}             \def\Gam{\Gamma}
\def\del{\delta}             \def\Del{\Delta}
\def\eps{\varepsilon}
\def\ep{\epsilon}
\def\zet{\zeta}
\def\tet{\theta}             \def\Tet{\Theta}
\def\iot{\iota}
\def\kap{\kappa}
\def\ka{\kappa}
\def\lam{\lambda}            \def\Lam{\Lambda}
\def\la{\lambda}             \def\La{\Lambda}
\def\sig{\sigma}             \def\Sig{\Sigma}
\def\si{\sigma}              \def\Si{\Sigma}
\def\vphi{\varphi}
\def\ome{\omega}             \def\Ome{\Omega}
\def\om{\omega}              \def\Om{\Omega}

\def\fA{{\mathfrak {A}}}
\def\fB{{\mathfrak {B}}}
\def\fC{{\mathfrak {C}}}
\def\fD{{\mathfrak {D}}}
\def\fE{{\mathfrak {E}}}
\def\fF{{\mathfrak {F}}}
\def\fG{{\mathfrak {G}}}
\def\fH{{\mathfrak {H}}}
\def\fI{{\mathfrak {I}}}
\def\fJ{{\mathfrak {J}}}
\def\fK{{\mathfrak {K}}}
\def\fL{{\mathfrak {L}}}
\def\fM{{\mathfrak {M}}}
\def\fN{{\mathfrak {N}}}
\def\fO{{\mathfrak {O}}}
\def\fP{{\mathfrak {P}}}
\def\fQ{{\mathfrak {Q}}}
\def\fR{{\mathfrak {R}}}
\def\fS{{\mathfrak {S}}}
\def\fT{{\mathfrak {T}}}
\def\fU{{\mathfrak {U}}}
\def\fV{{\mathfrak {V}}}
\def\fW{{\mathfrak {W}}}
\def\fX{{\mathfrak {X}}}
\def\fY{{\mathfrak {Y}}}
\def\fZ{{\mathfrak {Z}}}

\def\fa{{\mathfrak {a}}}
\def\fb{{\mathfrak {b}}}
\def\fc{{\mathfrak {c}}}
\def\fd{{\mathfrak {d}}}
\def\fe{{\mathfrak {e}}}
\def\ff{{\mathfrak {f}}}
\def\fg{{\mathfrak {g}}}
\def\fh{{\mathfrak {h}}}
%Attn:\fi can't be defined.
\def\fj{{\mathfrak {j}}}
\def\fk{{\mathfrak {k}}}
\def\fl{{\mathfrak {l}}}
\def\fm{{\mathfrak {m}}}
\def\fn{{\mathfrak {n}}}
\def\fo{{\mathfrak {o}}}
\def\fp{{\mathfrak {p}}}
\def\fq{{\mathfrak {q}}}
\def\fr{{\mathfrak {r}}}
\def\fs{{\mathfrak {s}}}
\def\ft{{\mathfrak {t}}}
\def\fu{{\mathfrak {u}}}
\def\fv{{\mathfrak {v}}}
\def\fw{{\mathfrak {w}}}
\def\fx{{\mathfrak {x}}}
\def\fy{{\mathfrak {y}}}
\def\fz{{\mathfrak {z}}}

\def\bbA{{\mathbb {A}}}
\def\bbB{{\mathbb {B}}}
\def\bbC{{\mathbb {C}}}
\def\bbD{{\mathbb {D}}}
\def\bbE{{\mathbb {E}}}
\def\bbF{{\mathbb {F}}}
\def\bbG{{\mathbb {G}}}
\def\bbH{{\mathbb {H}}}
\def\bbI{{\mathbb {I}}}
\def\bbJ{{\mathbb {J}}}
\def\bbK{{\mathbb {K}}}
\def\bbL{{\mathbb {L}}}
\def\bbM{{\mathbb {M}}}
\def\bbN{{\mathbb {N}}}
\def\bbO{{\mathbb {O}}}
\def\bbP{{\mathbb {P}}}
\def\bbQ{{\mathbb {Q}}}
\def\bbR{{\mathbb {R}}}
\def\bbS{{\mathbb {S}}}
\def\bbT{{\mathbb {T}}}
\def\bbU{{\mathbb {U}}}
\def\bbV{{\mathbb {V}}}
\def\bbW{{\mathbb {W}}}
\def\bbX{{\mathbb {X}}}
\def\bbY{{\mathbb {Y}}}
\def\bbZ{{\mathbb {Z}}}

\def\cA{{\mathcal {A}}}
\def\cB{{\mathcal {B}}}
\def\cC{{\mathcal {C}}}
\def\cD{{\mathcal {D}}}
\def\cE{{\mathcal {E}}}
\def\cF{{\mathcal {F}}}
\def\cG{{\mathcal {G}}}
\def\cH{{\mathcal {H}}}
\def\cI{{\mathcal {I}}}
\def\cJ{{\mathcal {J}}}
\def\cK{{\mathcal {K}}}
\def\cL{{\mathcal {L}}}
\def\cM{{\mathcal {M}}}
\def\cN{{\mathcal {N}}}
\def\cO{{\mathcal {O}}}
\def\cP{{\mathcal {P}}}
\def\cQ{{\mathcal {Q}}}
\def\cR{{\mathcal {R}}}
\def\cS{{\mathcal {S}}}
\def\cT{{\mathcal {T}}}
\def\cU{{\mathcal {U}}}
\def\cV{{\mathcal {V}}}
\def\cW{{\mathcal {W}}}
\def\cX{{\mathcal {X}}}
\def\cY{{\mathcal {Y}}}
\def\cZ{{\mathcal {Z}}}

\def\tA{{\widetilde{A}}}
\def\tB{{\widetilde{B}}}
\def\tC{{\widetilde{C}}}
\def\tD{{\widetilde{D}}}
\def\tE{{\widetilde{E}}}
\def\tF{{\widetilde{F}}}
\def\tG{{\widetilde{G}}}
\def\tH{{\widetilde{H}}}
\def\tI{{\widetilde{I}}}
\def\tJ{{\widetilde{J}}}
\def\tK{{\widetilde{K}}}
\def\tL{{\widetilde{L}}}
\def\tM{{\widetilde{M}}}
\def\tN{{\widetilde{N}}}
\def\tO{{\widetilde{O}}}
\def\tP{{\widetilde{P}}}
\def\tQ{{\widetilde{Q}}}
\def\tR{{\widetilde{R}}}
\def\tS{{\widetilde{S}}}
\def\tT{{\widetilde{T}}}
\def\tU{{\widetilde{U}}}
\def\tV{{\widetilde{V}}}
\def\tW{{\widetilde{W}}}
\def\tX{{\widetilde{X}}}
\def\tY{{\widetilde{Y}}}
\def\tZ{{\widetilde{Z}}}

\def\tcA{{\widetilde{\mathcal {A}}}}
\def\tcB{{\widetilde{\mathcal {B}}}}
\def\tcC{{\widetilde{\mathcal {C}}}}
\def\tcD{{\widetilde{\mathcal {D}}}}
\def\tcE{{\widetilde{\mathcal {E}}}}
\def\tcF{{\widetilde{\mathcal {F}}}}
\def\tcG{{\widetilde{\mathcal {G}}}}
\def\tcH{{\widetilde{\mathcal {H}}}}
\def\tcI{{\widetilde{\mathcal {I}}}}
\def\tcJ{{\widetilde{\mathcal {J}}}}
\def\tcK{{\widetilde{\mathcal {K}}}}
\def\tcL{{\widetilde{\mathcal {L}}}}
\def\tcM{{\widetilde{\mathcal {M}}}}
\def\tcN{{\widetilde{\mathcal {N}}}}
\def\tcO{{\widetilde{\mathcal {O}}}}
\def\tcP{{\widetilde{\mathcal {P}}}}
\def\tcQ{{\widetilde{\mathcal {Q}}}}
\def\tcR{{\widetilde{\mathcal {R}}}}
\def\tcS{{\widetilde{\mathcal {S}}}}
\def\tcT{{\widetilde{\mathcal {T}}}}
\def\tcU{{\widetilde{\mathcal {U}}}}
\def\tcV{{\widetilde{\mathcal {V}}}}
\def\tcW{{\widetilde{\mathcal {W}}}}
\def\tcX{{\widetilde{\mathcal {X}}}}
\def\tcY{{\widetilde{\mathcal {Y}}}}
\def\tcZ{{\widetilde{\mathcal {Z}}}}

\def\tfA{{\widetilde{\mathfrak {A}}}}
\def\tfB{{\widetilde{\mathfrak {B}}}}
\def\tfC{{\widetilde{\mathfrak {C}}}}
\def\tfD{{\widetilde{\mathfrak {D}}}}
\def\tfE{{\widetilde{\mathfrak {E}}}}
\def\tfF{{\widetilde{\mathfrak {F}}}}
\def\tfG{{\widetilde{\mathfrak {G}}}}
\def\tfH{{\widetilde{\mathfrak {H}}}}
\def\tfI{{\widetilde{\mathfrak {I}}}}
\def\tfJ{{\widetilde{\mathfrak {J}}}}
\def\tfK{{\widetilde{\mathfrak {K}}}}
\def\tfL{{\widetilde{\mathfrak {L}}}}
\def\tfM{{\widetilde{\mathfrak {M}}}}
\def\tfN{{\widetilde{\mathfrak {N}}}}
\def\tfO{{\widetilde{\mathfrak {O}}}}
\def\tfP{{\widetilde{\mathfrak {P}}}}
\def\tfQ{{\widetilde{\mathfrak {Q}}}}
\def\tfR{{\widetilde{\mathfrak {R}}}}
\def\tfS{{\widetilde{\mathfrak {S}}}}
\def\tfT{{\widetilde{\mathfrak {T}}}}
\def\tfU{{\widetilde{\mathfrak {U}}}}
\def\tfV{{\widetilde{\mathfrak {V}}}}
\def\tfW{{\widetilde{\mathfrak {W}}}}
\def\tfX{{\widetilde{\mathfrak {X}}}}
\def\tfY{{\widetilde{\mathfrak {Y}}}}
\def\tfZ{{\widetilde{\mathfrak {Z}}}}

%roman letters with a tilde
\def\Atil{{\widetilde A}}          \def\atil{{\tilde a}}
\def\Btil{{\widetilde B}}          \def\btil{{\tilde b}}
\def\Ctil{{\widetilde C}}          \def\ctil{{\tilde c}}
\def\Dtil{{\widetilde D}}          \def\dtil{{\tilde d}}
\def\Etil{{\widetilde E}}          \def\etil{{\tilde e}}
\def\Ftil{{\widetilde F}}          \def\ftil{{\tilde f}}
\def\Gtil{{\widetilde G}}          \def\gtil{{\tilde g}}
\def\Htil{{\widetilde H}}          \def\htil{{\tilde h}}
\def\Itil{{\widetilde I}}          \def\itil{{\tilde i}}
\def\Jtil{{\widetilde J}}          \def\jtil{{\tilde j}}
\def\Ktil{{\widetilde K}}          \def\ktil{{\tilde k}}
\def\Ltil{{\widetilde L}}          \def\ltil{{\tilde l}}
\def\Mtil{{\widetilde M}}          \def\mtil{{\tilde m}}
\def\Ntil{{\widetilde N}}          \def\ntil{{\tilde n}}
\def\Otil{{\widetilde O}}          \def\otil{{\tilde o}}
\def\Ptil{{\widetilde P}}          \def\ptil{{\tilde p}}
\def\Qtil{{\widetilde Q}}          \def\qtil{{\tilde q}}
\def\Rtil{{\widetilde R}}          \def\rtil{{\tilde r}}
\def\Stil{{\widetilde S}}          \def\stil{{\tilde s}}
\def\Ttil{{\widetilde T}}          \def\ttil{{\tilde t}}
\def\Util{{\widetilde U}}          \def\util{{\tilde u}}
\def\Vtil{{\widetilde V}}          \def\vtil{{\tilde v}}
\def\Wtil{{\widetilde W}}          \def\wtil{{\tilde w}}
\def\Xtil{{\widetilde X}}          \def\xtil{{\tilde x}}
\def\Ytil{{\widetilde Y}}          \def\ytil{{\tilde y}}
\def\Ztil{{\widetilde Z}}          \def\ztil{{\tilde z}}

%roman letters with a bar
%\def\abar{{\bar a}} \def\Abar{{\bar A}}
%\def\bbar{{\bar b}} \def\Bbar{{\bar B}}
%\def\cbar{{\bar c}} \def\Cbar{{\bar C}}
%\def\dbar{{\bar d}} \def\Dbar{{\bar D}}
%\def\ebar{{\bar e}} \def\Ebar{{\bar E}}
%\def\fbar{{\bar f}} \def\Fbar{{\bar F}}
%\def\gbar{{\bar g}} \def\Gbar{{\bar G}}
%\def\hBar{{\bar h}} \def\Hbar{{\bar H}}
%\def\ibar{{\bar i}} \def\Ibar{{\bar I}}
%\def\jbar{{\bar j}} \def\Jbar{{\bar J}}
%\def\kbar{{\bar k}} \def\Kbar{{\bar K}}
%\def\lbar{{\bar l}} \def\Lbar{{\bar L}}
%\def\mbar{{\bar m}} \def\Mbar{{\bar M}}
%\def\nbar{{\bar n}} \def\Nbar{{\bar N}}
%\def\obar{{\bar o}} \def\Obar{{\bar O}}
%\def\pbar{{\bar p}} \def\Pbar{{\bar P}}
%\def\qbar{{\bar q}} \def\Qbar{{\bar Q}}
%\def\rbar{{\bar r}} \def\Rbar{{\bar R}}
%\def\sbar{{\bar s}} \def\Sbar{{\bar S}}
%\def\tbar{{\bar t}} \def\Tbar{{\bar T}}
%\def\ubar{{\bar u}} \def\Ubar{{\bar U}}
%\def\vbar{{\bar v}} \def\Vbar{{\bar V}}
%\def\wbar{{\bar w}} \def\Wbar{{\bar W}}
%\def\xbar{{\bar x}} \def\Xbar{{\bar X}}
%\def\ybar{{\bar y}} \def\Ybar{{\bar Y}}
%\def\zbar{{\bar z}}  \def\Zbar{{\bar Z}}

%roman letters with a hat
\def\ahat{{\hat a}}          \def\Ahat{{\widehat A}}
\def\bhat{{\hat b}}          \def\Bhat{{\widehat B}}
\def\chat{{\hat c}}          \def\Chat{{\widehat C}}
\def\dhat{{\hat d}}          \def\Dhat{{\widehat D}}
\def\ehat{{\hat e}}          \def\Ehat{{\widehat E}}
\def\fhat{{\hat f}}          \def\Fhat{{\widehat F}}
\def\ghat{{\hat g}}          \def\Ghat{{\widehat G}}
\def\hhat{{\hat h}}          \def\Hhat{{\widehat H}}
\def\ihat{{\hat i}}          \def\Ihat{{\widehat I}}
\def\jhat{{\hat j}}          \def\Jhat{{\widehat J}}
\def\khat{{\hat k}}          \def\Khat{{\widehat K}}
\def\lhat{{\hat l}}          \def\Lhat{{\widehat L}}
\def\mhat{{\hat m}}          \def\Mhat{{\widehat M}}
\def\nhat{{\hat n}}          \def\Nhat{{\widehat N}}
\def\ohat{{\hat o}}          \def\Ohat{{\widehat O}}
\def\phat{{\hat p}}          \def\Phat{{\widehat P}}
\def\qhat{{\hat q}}          \def\Qhat{{\widehat Q}}
\def\rhat{{\hat r}}          \def\Rhat{{\widehat R}}
\def\shat{{\hat s}}          \def\Shat{{\widehat S}}
\def\that{{\hat t}}          \def\That{{\widehat T}}
\def\uhat{{\hat u}}          \def\Uhat{{\widehat U}}
\def\vhat{{\hat v}}          \def\Vhat{{\widehat V}}
\def\what{{\hat w}}          \def\What{{\widehat W}}
\def\xhat{{\hat x}}          \def\Xhat{{\widehat X}}
\def\yhat{{\hat y}}          \def\Yhat{{\widehat Y}}
\def\zhat{{\hat z}}          \def\Zhat{{\widehat Z}}

\def\tg{{\widetilde g}}

\author{Michael Christ, Loukas Grafakos,  Petr Honz\'\i k,
 Andreas Seeger}
 \address{Michael Christ\\ Department of Mathematics
\\ University of California\\
Berkeley, CA 94720-3840, USA} \email{mchrist@math.berkeley.edu}

\address{
        Loukas Grafakos\\
        Department of Mathematics\\
        University of Missouri \\
        Co- lumbia,  MO 65211, USA}
\email{loukas@math.missouri.edu}

\address{Petr Honz\'\i k\\
        Department of Mathematics\\
        University of Missouri \\
        Columbia, MO 65211, USA}
\email{honzikp@math.missouri.edu}

\address{
        Andreas Seeger\\
        Department of Mathematics\\
        University of Wisconsin \\
        Madison, WI 53706, USA}
\email{seeger@math.wisc.edu}
\thanks{Christ, Grafakos and Seeger 
were supported in part by  NSF grants.
Honz\'\i k was supported by 201/03/0931 Grant Agency of the
Czech Republic}

%\date{\today}
\date{March 7, 2004}

\title[Maximal functions and Mikhlin-H\"ormander multipliers]
% runningheader
{Maximal functions associated
with Fourier multipliers of Mikhlin-H\"ormander type}

  % Title

\begin{abstract}
We  show that  maximal operators
formed by dilations of   Mikhlin-H\"ormander multipliers
% transformations of
are typically not bounded on $L^p(\bbR^d)$. We also give rather
weak  conditions in terms of the decay
 of such   multipliers under which $L^p$ boundedness of the maximal operators
holds.
\end{abstract}

\maketitle

\section{Introduction}\label{intro}

%%\centerline{\bf Introduction}
%%\medskip

For a bounded  Fourier multiplier $m$
on $\mathbb R^d$ and a Schwartz function
$f$ in $ \mathcal S(\mathbb R^d)$ define  the
  maximal function associated with $m$
by
$$\mathcal M_m f(x)  = \sup_{t>0}
\big|\mathcal F^{-1}[m(t\cdot) \widehat f \,](x) \big|.$$

We are interested in the class of multipliers that satisfy the
estimates of the standard Mikhlin-H\"ormander multiplier theorem
\begin{equation}\label{1}
|\partial^\alpha m(\xi)|\leq C_\alpha |\xi|^{-\alpha}
\end{equation}
for all  (or  sufficiently large)
 multiindices $\alpha$. More precisely,
let $L^r_\gamma$ be the standard Bessel-potential (or Sobolev) space
with norm
$$\|f\|_{L^r_\gamma}=\|(I-\Delta)^{\gamma/2} f \|_r;$$
here we include the case $r=1$. Let
$\phi$ be a smooth function supported in
$\{\xi: 1/2< |\xi|< 2\}$ which is nonvanishing on
$\{\xi: 1/\sqrt 2\le |\xi|\le\sqrt 2\}$.
% $L^2_\gamma$ is the standard $L^2$ Sobolev space.
Then one
 imposes conditions on $m$ of the form
\begin{equation}\label{2}
\sup_{k\in \mathbb Z} \|\phi m(2^k \cdot)\|_{L^r_\gamma}<\infty.
\end{equation}

The function $m$ is a Fourier multiplier on all $L^p$,
$1<p<\infty$ if \eqref{2} holds for $\gamma>d/r$, with $1\le r\le 2$  and
the condition for $r=2$ is the least restrictive one
(see \cite{Hoe}). Concerning the maximal operator
Dappa and Trebels \cite{DT} showed using Calder\'on-Zygmund theory that if
$\cM_m$ is a priori bounded on some $L^q$, $q>1$ and if \eqref{2} holds
for $r=1$, $\gamma>d,$ then
$\cM_m$ is of weak type $(1,1)$ and thus bounded on $L^p$ for $1<p<q$.
Using square function   estimates, the
 $L^2$ boundedness
of $\mathcal M_m$ has been shown  in    \cite{C},  \cite{DT}
under certain additional decay assumptions
  ({\it cf.} also  \cite{RdF}). For instance,  it follows from
 \cite{DT} that
\begin{equation}\label{3}
\|\mathcal M_m f\|_p
\le C_{p} \Big(\sum_{k\in \mathbb Z}
\|\phi m(2^k \cdot)\|_X^2
%{L^1_{\gamma}}^2
\Big)^{1/2}
%{L^2_\gamma}^2\Big)^{1/2}
 \|f\|_p,
\end{equation}
 with $X=L^p_{d/p+\epsilon}$ for $1<p\le 2$ 
, and with $X=L^2_{d/2+\epsilon}$ ,  for $2\le p<\infty$. Further results in terms
of weaker  differentiability assumptions are in \cite {C},
\cite{DT}, especially for classes of radial multipliers. Moreover,
if $m$
%is smooth in $\bbR^{d}\setminus 0$ and
is homogeneous of
degree $0$ then trivially $|\cM_m f|= |\cF^{-1}[m\widehat f \,]|$;
 this observation can be used to build more general
classes of symbols without decay assumptions for which $\cM_m$ is
$L^p$ bounded.

A problem left open in \cite{DT} is
%whether this type of decay assumption on $m$
%is necessary or
whether the Mikhlin-H\"ormander type assumption  in
\eqref{1} or    \eqref{2}  alone is
sufficient to prove boundedness of the maximal operator $\mathcal M_m$.
We show here that
some additional  assumption is needed; indeed this applies
already  to the
dyadic  maximal function associated with $m$,
defined by
\begin{equation}
\label{dyadic-maximal}
M_m f  = \sup_{k\in \mathbb Z}
\big|\mathcal F^{-1}[m(2^k\cdot) \widehat f \,]\big|,
\end{equation}
which of course is dominated by  $\mathcal M_mf $.
\medskip

\noindent{\bf Example.}
\emph{Let
 $\{v(l)\}_{l=0}^\infty$ be a
positive increasing and unbounded  sequence.
%with $\lim_{k\to \infty} v(k)=\infty$.
 Then there is a Fourier multiplier $m$ satisfying
\begin{equation}\label{example}
\sup_\xi \big|\partial_\xi^\alpha\big(\phi(\xi) m(2^k
\xi)\big)\big| \le C_\alpha \frac{ v(|k|)}{\sqrt{ \log(|k|+2)}},
\quad k\in \bbZ,
\end{equation}
with $C_\alpha<\infty$ for all  multiindices $\alpha$,
so that the associated dyadic maximal operator
$M_m$  is unbounded on $L^p(\mathbb R^d)$
 for  $1<p<\infty$.
}

\medskip

This counterexample will be explicitly constructed in $\S2$.
Taking $v(l)= \sqrt{ \log(l+2)}$ we see
 that  there exists $m$ satisfying  \eqref{1},
so that $M_m$, and hence $\mathcal M_m$, are
unbounded on   $L^p(\mathbb R^d)$ for $1<p<\infty$.
In view of these examples it is not unexpected that
unboundedness of $M_m$ holds in fact  for the {\it typical} multiplier
 satisfying  \eqref{1}, {\it i.e.} on a residual set in the
sense of Baire category.
In order to formulate a result let
 $\fS$ be the  space of
functions $m\in C^\infty(\bbR^d\setminus \{0\})$ satisfying \eqref{1}
with $C_\alpha<\infty$ for all  multiindices $\alpha$.
It is easy to see that  $\fS$ is a
Fr\'echet-space with the  topology given by the countable  family of norms
\begin{equation}
\label{seminorms}
\|m\|_{(j)}= \sup_{|\alpha|\le j} \sup_{\xi\in \bbR^d}|\xi|^{|\alpha|}
|\partial_\xi^\alpha m(\xi)|.
\end{equation}

%\begin{equation}
%\big|\partial_\xi^\alpha\big(\phi(\xi) m(2^k \xi)\big)\big| \le C_\alpha
%\frac{
%v(|k|)}{\sqrt{
%\log(|k|+2)}}
%\end{equation}

%The best constants in \eqref{example} define seminorms, and it is not hard to
% see that $\fS(w)$ becomes a Fr\'echet space with this choice of seminorms.

Let $\mathcal S_0$ denote the space of
 Schwartz functions whose Fourier transform have compact
support in $\bbR^d\setminus\{0\}$ and let
$\fS^M$ be  the space of  all $m\in \fS$ for which
$$
\sup\{\|M_m f\|_p: f\in \mathcal S_0, \|f\|_p\le 1\}
$$
is finite for some $p\in (1,\infty)$. Thus  $m\in \fS^M$
if and only if the linear operator $f\mapsto\{\cF^{-1}[m(t\cdot)\widehat f\,]\}_{t>0}$ extends to a bounded operator from $L^p$ to $L^p(L^\infty)$ for some $p\in (1,\infty)$; in other words $m\in \fS^M$ if and only if
$M_m$ extends to a bounded operator on $L^p(\bbR^d)$, for some
$p\in (1,\infty)$.
\medskip

\noindent{\bf Theorem 1.1.}
\emph{$\fS^{M}$ is of  first category in $\fS$, in the sense of Baire.}

\medskip

In terms of positive results we note that there is a
 significant gap bet\-ween the known   conditions in \eqref{3} and
the weak decay \eqref{example}.
Assuming $\|\phi m(2^k \cdot)\|_{L^1_{d+\epsilon}}= O(|k|^{-\alpha})$,
then  \eqref{3} yields
$L^p$ boundedness for  $1<p <\infty$ only when
$\alpha>1/2$.
 We shall see that
 this result
remains in fact  valid
% for all $\alpha>0$ and
   under the weaker assumption
\begin{equation}
\label{logpowerone} \|\phi m(2^k \cdot)\|_{L^1_{d+\epsilon}}
\lc (\log (|k|+2))^{-1-\epsilon}.
\end{equation} In what follows we shall
 mainly aim for minimal decay but will also try to formulate reasonable 
smoothness assumptions.

To formulate a general result we recall the definition of
 the nonincreasing rearrangement of a sequence $\omega$,
defined for $t\ge 0$ by
$$\omega^*(t)= \sup\big \{\la>0:
\card\big(\{k: |\omega(k)|>\lambda\}\big)> t\big\};
$$
note that $\omega^*(0)=\sup_k|\omega(k)|$ and $\omega^* $
is constant on the intervals $[n,n+1)$, $n=0,1,2,\dots$.
\medskip

\noindent{\bf Theorem 1.2.}
\emph{Let $1<p<\infty$, $1/p+1/p'=1$ and let $\omega:\bbZ\to [0,\infty)$ 
satisfy}
\begin{equation}
\label{omega}
\omega^*(0)+\sum_{l=1}^\infty\frac{\om^*(l)}{l}  <\infty.
\end{equation}

\emph{ (i) Suppose that for some 
$\alpha>d/p$ we have}
\begin{equation}
\label{thm12assumption}
 \Big(\int_{\bbR^d} \big|\cF^{-1}[\phi m(2^k\cdot)]\big|^{p'}(1+|x|)^{\alpha p'} 
dx\Big)^{1/p'} \le \omega(k), \quad k\in\bbZ \, ,
\end{equation}
\emph{then $\cM_m$ is bounded on $L^p(\bbR^d)$.}

\emph{ (ii) If \eqref{thm12assumption} holds for $p'=1$ then 
$\cM_m$ maps $L^\infty$ to $BMO$.}

\emph{ (iii) If for some $\eps>0$ }
\begin{equation}
\label{thm12assumptionweaktype}
\sup_x (1+|x|)^{d+\epsilon} \big| \cF^{-1}[\phi m(2^k\cdot)](x)\big |
\le \omega(k), \quad k\in\bbZ.
\end{equation}
\emph{Then  $\cM_m$ is of weak type $(1,1)$, and $\cM_m$ maps $H^1$ to $L^1$.}

\medskip

By the  Hausdorff-Young inequality for $p\le 2$ one  deduces

\medskip 

\noindent{\bf Corollary 1.3.} \emph{Suppose $1<p<\infty$, 
 $r=\min \{p,2\}$, and $\alpha>d/r$. Suppose that}
\begin{equation}
\label{sob-Lr}
 \|\phi m(2^k \cdot)\|_{L^r_\alpha}\le \omega(k), \quad k\in \bbZ,
\end{equation}
\emph{where  $\omega$ satisfies \eqref{omega}. Then $\cM_m$ 
is bounded on $L^p(\bbR^d)$.}

%\emph{ If these assumptions hold for $r=1$ 
%then  $\cM_m$ is of weak type $(1,1)$.}

\medskip

%We note that by the results in \cite{DT} the condition
%$\|\phi m(2^k \cdot)\|_{L^1_\gamma}=O(1)$ and the a
%priori boundedness of $\cM_m$ on some $L^q$ implies the
%weak type (1,1) property and
%thus $L^p$ boundedness for $1<p<q$.
%Thus we have to show only  $L^p$ boundedness for a range $p\in [q,\infty)$.

In particular
we conclude that the condition
\begin{equation}
\label{lq}
\Big(\sum_{k\in \bbZ} \|\phi m(2^k\cdot)\|_{L^r_{\alpha}}^q\Big)^{1/q}<\infty
\end{equation}
with $\alpha$, $r$ as in the corollary, implies $L^p$ boundedness. Indeed 
\eqref{lq} implies that
$\omega^*(l)=O(l^{-1/q})$ as $l\to \infty.$ Of course $L^p$ boundedness  
also holds if $\omega^*(l)\lc (\log(2+l))^{-1-\eps}$ etc.
which covers
condition \eqref{logpowerone}.

\medskip

Finally we state  a more elementary but closely related result
about maximal functions for a finite number of H\"ormander-Mikhlin type
multipliers $m_\nu$, with no decay assumptions and not necessarily
generated by dilating a single multiplier.

\medskip

\noindent{\bf Theorem 1.4.}
\emph{Let $1<p<\infty$, $1/p+1/p'=1$
%$1/r= 1/2+1/2p$,
%%$r=\min\{p,2\}$ and $\alpha>d/r$, 
and let $\{m_\nu\}_{\nu\ge 1}$ be  a sequence of multipliers and define a maximal operator by}
 \begin{equation*}
\fM_n f(x)  =\sup_{1\le \nu\le n}\big |\cF^{-1}[m_\nu \widehat f](x)\big|.
\end{equation*}
\emph{Suppose that}
\begin{equation}
\label{thm14assumption}
%%\label{sob-Lr-nu}
%\sup_{ \nu}\sup_{k\in \bbZ}
%  \|\phi m_\nu(2^k \cdot)\|_{L^r_\alpha}\le A.
%\end{equation}
\sup_\nu \sup_{k\in \bbZ}
\Big(\int_{\bbR^d}  \big|\cF^{-1}[\phi m_\nu(2^k\cdot)]\big|^{p'}(1+|x|)^{\alpha p'} 
dx\Big)^{1/p'}\le A
\end{equation}
for some 
$\alpha>d/p$.
\emph{
Then for $f\in L^p(\bbR^d)$}
\begin{equation}\label{lognbound}
\|\fM_n f \|_p\le C_p A \log (n+1)\|f\|_p.
\end{equation}

%\emph{If  for some $\eps>0$}
%\begin{equation}\label{thm14assumptionweaktype}
%\sup_\nu \sup_{k\in \bbZ}\sup_x (1+|x|)^{d+\eps}
%\big|\cF^{-1}[\phi m_\nu(2^k\cdot)]\big|\le A
%\end{equation}
%\emph{then there is the weak type inequality}
%$$
%\sup_{\la>0}\la \,\meas\big(\big\{x:\fM_n f(x)>\la\big\}\big)\le C_\eps 
%\log(n+1) \|f\|_1.
%$$
%\medskip

Again if the above assumptions hold for $p=1$
then  a weak type (1,1) inequality and an $H^1\to L^1$ inequality 
 hold, and if $p=\infty$ we have an $L^\infty\to BMO$ inequality, all with
 constant $O(\log (n+1))$.

\medskip

{\it Structure of the paper.} In \S2 we shall provide the above
mentioned examples for unboundedness and prove Theorem 1.1. A
tiling lemma for finite sets of integers  and other preliminaries
needed in the proof of Theorem 1.2 are provided in \S3. \S4
contains the main relevant estimates for multipliers supported in
a finite union of annuli. In \S5 we conclude the proof of Theorem
1.2 and in \S6 we give the proof of Theorem 1.4. Finally we state
some extensions and open problems.

\section{Unboundedness of the maximal operator}
\label{counterexample} We shall explicitly construct an example
satisfying \eqref{example} and then use our example to prove
Theorem 1.1.

Define
$
S=\{1,-1,i,-i\}
$
and let $S^N$ be  the set of sequences of length $N$ on $S$.
Enumerate the  $4^N$ elements in $S^N$ by
 $\{s_\ka\}_{\ka=1}^{4^N}$.
Let  $\Phi$ be a smooth function supported in $ 3/4\le |\xi| \le 5/4 $, so that $\Phi(\xi)=1$
whenever $   7/8\le |\xi| \le 9/8 $.
We let
$$
m_N(\xi):=
\sum_{\kappa=1}^{4^N}\sum_{j=1}^{N} s_{\ka}(j)\Phi(2^{-N\ka-j}\xi)
$$
which is supported in $\{\xi:1/2\le |\xi|\le 2^{N4^{N}+N+1}\}$,
and define $m$ by
  \begin{equation}
\label{counterex}
m(\xi) =\sum_{N=1}^\infty N^{-1/2} v(4^N)
m_N(2^{-N8^{N}} \xi).
\end{equation}
One observes that the terms in this sum have disjoint supports and
that  $m$ satisfies
condition \eqref{example}.

Fix $1<p<\infty$.
We will  test the maximal operator $ M_{m}$ on
functions  $ f_{N,p} $ defined as follows.
Pick a Schwartz function $\Psi$ such that $\|\Psi\|_p=1$
and so that $\supp \widehat \Psi $ is contained in the ball $|\xi|\le 1/8$.
 For  $1\le j\le N$
define
$$
g_N(x)=\sum_{j=1}^N e^{2\pi i 2^j x_1}\Psi(x),
$$
and set
\begin{equation}\label{fN}
f_{N,p}(x) = N^{-1/2}2^{dN8^N/p} g_N(2^{N8^N} x).
\end{equation}
Then $\widehat {g_N} (\xi )=  \sum_{j=1}^N \widehat \Psi(\xi-2^j e_1)$
% and by orthogonality  $\|g_N\|_2=N^{1/2}$.
and, by   Littlewood-Paley theory,
$$\|g_N\|_p\le c_p N^{1/2}, \quad 1<p<\infty.$$
Thus
$$\|f_{N,p}\|_p\le c_p <\infty, \quad 1<p<\infty,$$
uniformly in $N$.

The main observation is
\begin{equation} \label{5}
\Big\|\sup_{1\le k\le N4^N}
\Big| \mathcal F^{-1}[ m_N(2^k \cdot ) \widehat {g_N}]\Big|\,
\Big\|_p\ge C N.
\end{equation}

Given \eqref{5} we quickly derive the asserted unboundedness of $M_m$.
Namely, by the support properties of the $m_n$ it follows that
$$
m_n(2^{k-n8^n}\xi) \widehat {g_N}(2^{-N8^N}\xi)
=0 \quad\text{ if } N\neq n, \quad 1\le k\le N4^N.
$$
Thus, setting  $a_n= n^{-1/2} v(4^n) $, we obtain
\begin{align*}
M_m f_{N,p}(x)&\ge \sup_{1\le k\le N4^N}\Big|\sum_{n=1}^\infty a_n \mathcal F^{-1}[  m_n(2^{k-n 8^n}\cdot) \widehat {f_{N,p}}](x)\Big|
\\
&\ge \sup_{1\le k\le N4^N} a_N N^{-\frac12}
\big | \mathcal F^{-1}[  m_N(2^{k-N 8^N}\cdot ) 2^{-\frac{dN 8^N}{p'}} \widehat {g_N}
(2^{-N 8^N }\cdot)](x)\big |
\\
&=
a_N N^{-\frac12}2^{\frac{dN 8^N}{p}}
\sup_{1\le k\le N4^N}\big|
\mathcal F^{-1}[m_N (2^k\cdot)  \widehat{ g_N}]
(2^{N8^N}x) \big|.
\end{align*}
Taking $L^p(\mathbb R^d)$ norms and using \eqref{5} we conclude that
\begin{equation}
\label{counterex-conclusion}
\|M_m f_{N,p}\|_p\ge  C a_N N^{1/2} = Cv(4^N).
\end{equation}
By the assumed
unboundedness of the increasing sequence  $v$ it follows that $M_m$
is not bounded on $L^p$.

\medskip

\emph{ Proof of \eqref{5}.}
For any complex number $z$ the quantity
 $\sup_{c\in S} {\rm Re}\; (cz)$ is at least
$ |z|/\sqrt{ 2}$.
Thus for  $x\in \mathbb R^d$  and $1\leq j\leq N$ we may pick $c_j(x)\in S$ such that
\begin{equation} \label{6}
{\rm Re}\; \big(c_j(x) e^{2\pi i 2^j x_1}\Psi(x) \big)\geq  |\Psi(x)|/\sqrt {2 }.
\end{equation}
We can find     $\kappa_x$ in $\{1,\dots , 4^{N}\}$ such that
%$\{c_j(x)\}_{j=1}^{N}=s_{\kappa_x}$.
$$
c_j(x)=s_{\kappa_x}(j), \quad j=1,\dots, N
.$$
Taking $k=\kappa_x N$ we obtain
\begin{align*}
\sup_{1\le k\le N4^N} &\big|\mathcal F^{-1}[ m_N(2^k \cdot ) \widehat {g_N}]\big(x)|\\&
\geq
{\rm Re} \int \sum_{l=1}^{4^N}\sum_{\nu=1}^{N}s_l(\nu)
\Phi(2^{-Nl-\nu}2^{N\kappa_x}\xi)
\sum_{j=1}^N \widehat \Psi(\xi-2^je_1 ) e^{2\pi i \inn{x}{\xi}}d\xi.
\end{align*}
Since $1\leq j \leq N$, the supports of $\Phi(2^{-Nl-\nu}2^{N\ka_x}\xi)$
and  $\widehat\Psi(\xi- 2^je_1)$ intersect only when $l=\kappa_x$ and
$j=\nu$.
 In this  case
$\Phi(2^{-Nl-\nu}2^{N\kappa_x}\xi)=\Phi(2^{-j} \xi)$ is equal to $1$
on the support of $\widehat{\Psi}
(\xi-2^j)$.
Therefore  we   obtain from \eqref{6} the pointwise estimate
\begin{align*}
\sup_{1\le k\le N4^N} &\big|\mathcal F^{-1}[ m_N(2^k \cdot )
 \widehat {g_N}]\big(x)|
\\&\geq
\sum_{j=1}^N{\rm Re }\,  \Big(s_{\ka_x}(j) \mathcal F^{-1}[ \widehat \Psi(\cdot-2^j)](x)\Big)
\geq  N |\Psi(x)|/\sqrt{2}\, .
\end{align*}
Taking $L^p$ norms yields  \eqref{5}. \qed

\medskip

\noi{\bf Proof of Theorem 1.1.}
The space $\fS $ is  a complete metric space and
 the metric is given by
$$d(m_1,m_2)=\sum_{j=0}^\infty 2^{-j}
\frac{\|m_1-m_2\|_{(j)}}{1+\|m_1-m_2\|_{(j)}}
$$ where $\|\cdot \|_{(j)}$ is defined in \eqref{seminorms}.

Let $f_{N,r'}$ be as in \eqref{fN} (with $p=r'$)  and for
 integers $r, n, N$,   all $\ge 2$,
consider the set
$$
\fS(r,n, N)=\{m\in \fS: \|M_{m} f_{N,r'}\|_{r'}\le n\},
$$
here $r'=r/(r-1)$,
and the set $$\fS(r,n)=\bigcap_{N=2}^\infty \fS(r,n,N).$$
We shall show that $\fS(r,n)$ is closed in $\fS$, and nowhere dense.
We also observe that
\begin{equation}
\label{fSM}
\fS^M\, \subset\, \bigcup_{r=2}^\infty \bigcup_{n=2}^\infty
\fS(r,n);
\end{equation}
thus $\fS^M$ is of first category. To see  \eqref {fSM} assume that
$M_m$ is bounded on $L^{p_0}$, for some $p_0>1$.
By the  theorem by Dappa and Trebels
mentioned before ({\it cf.} Proposition 3.2 below)
it follows that $M_m$ is bounded on
$L^p$ for  $1<p<p_0$, in particular bounded on $L^{r'}$
for some integer $r\ge 2$.
We note that  $f_{N,r'}\in \cS_0$ is such that
$\|f_{N,r'}\|_{r'}\le C_r$, independently of $N$.
Thus $m\in \fS(r,n)$ for sufficiently large $n$.

Next, in order to show that the sets $\fS(r,n)$ are  closed it suffices to show that the sets
$\fS(r,n,N)$ are closed for all $N\ge 2$.
For integers $l_1\le l_2$ denote by $\cS(l_1,l_2)$ the class of
Schwartz functions whose Fourier transform is supported in the annulus
$\{\xi: 2^{l_1-1}\le|\xi|\le 2^{l_2+1}\}$.
We observe the following inequality
$$
\|M_m f\|_{p}\le C(p) \|m\|_{(d+1)} (1+|l_2-l_1|) \|f\|_p, \quad \text{ if }f\in \cS(l_1,l_2),
$$
which (in view of the dependence on $l_1$, $l_2$) can be obtained by
standard techniques, see e.g. \cite{DT} or \cite{RdF}.
Note that every $f_{N,r'}$ is in some class
$\cS(l_1,l_2)$ with $l_2-l_1\le N$.
Now, if $m_\nu\in \fS(r, n,N)$ and $\lim_{\nu\to \infty} d(m_\nu, m)=0$ then
$$
\|M_m f_{N,r'}\|_{r'}\le  n+
\|M_{m-m_\nu} f_{N,r'}\|_{r'} \le n+C(r') \|m-m_\nu\|_{(d+1)} \|f_{N,r'}\|_{r'}
$$
 and since
$\|m-m_\nu\|_{(d+1)} \to 0$ we see that $m\in\fS(r,n,N)$.
Thus
$\fS(r,n,N)$ is closed.

Finally we need to  show that
 $\fS(r,n)$ is nowhere dense in  $\fS$; since this set is closed
we need to show that it  does not contain any open
balls.
Now if  $g\in \fS(r,n)$ then consider  the sequence $g_\nu=g+2^{-\nu} m$
where $m$ is as in \eqref{counterex}. Clearly $d(g_\nu, g)\to 0$.
However
by \eqref{counterex-conclusion} we have that
$g_\nu\notin \fS(r,n, N)$ for sufficiently large $N$ and thus
 $g_\nu\notin \fS(r,n)$ for any $r$, $n$.
Thus $\fS(r,n) $ is nowhere dense.\qed

\section{Preliminaries.}

\noi{\bf A tiling lemma.} In \S5 below we shall decompose the multiplier
into pieces with compact
but large support.
In order to  effectively estimate the maximal function associated
to these pieces
we shall use  the following
``tiling'' lemma for integers.

\medskip

\noindent{\bf Lemma 3.1.}
\emph{Let $N>0$ and let $E$ be a set of
%nonnegative
integers
of cardinality $\leq 2^{N}$.
Then we can find a set  $B=\{b_i\}_{i\in \mathbb Z}$ of integers,
such that}

\emph{
(i) the sets $b_i+E$ are pairwise disjoint,}

\emph{(ii) $b_i \in [i4^{N+1}, (i+1) 4^{N+1})$, and}

\emph {(iii) $\mathbb Z=\cup_{n=-4^{N+1}}^{4^{N+1}} (n+B)$.}

\medskip

\emph{Proof:}  Clearly (iii) is an immediate consequence of (ii).
We enumerate the set $E=\{e_\nu\}_{\nu=1}^{2^{N}}$.

We set $b_0=0$, and construct $b_j$, $b_{-j}$  for $j>0$ by induction.
Assume that
$b_i \in [i4^{N+1}, (i+1) 4^{N+1})$
has been constructed for $-j<i<j$ so that the sets
$b_i+E$ are pairwise disjoint.

For $\nu=1,\dots, 2^{N}$ we denote by $C_\nu^j$
the subset of all integers $c$ in
$ [j4^{N+1},
(j+1) 4^{N+1})$ with the property that
$e_\nu+c\in\cup_{i=1-j}^{j-1}(b_i+E).$

We shall verify
\begin{equation}\label{card}
\card\big(\cup_{\nu=1}^{2^{N}}C_\nu^j\big) \le 2^{2N+1}<4^{N+1}.
\end{equation}

Given \eqref{card} we may  simply take
$$b_j \in
[j4^{N+1}, (j+1) 4^{N+1}) \setminus  \big(\cup_{\nu=1}^{2^{N}}C_\nu^j\big)
$$
and by construction the sets $b_{1-j}+E, \dots, b_j+E$ are disjoint.

In order to verify \eqref{card} observe that
$e_\nu+c\in [j4^{N+1}+e_\nu,
(j+1) 4^{N+1}+e_\nu )$ if $c\in C_\nu^j$.
Thus $$\card(C_\nu^j)=
\card\big([j4^{N+1}+e_\nu,
(j+1) 4^{N+1}+e_\nu\big)\cap \cup_{i=1-j}^{j-1}(b_i+E)\big).$$
Since  by the induction assumption
$b_{i+2}-b_i>4^{N+1}$, if $i\ge 1-j$, $i+2\le j-1$, this  gives
$$\card\big(
[j4^{N+1}+e_\nu,
(j+1) 4^{N+1}+e_\nu)\cap \cup_{i=1-j}^{j-1}(b_i+\{e_\nu\})\big) \leq 2$$
for all $\nu$. This means
$\card(C_\nu^j)\leq 2\card(E)\leq 2^{N+1}$ and thus the cardinality of
$\cup_{\nu=1}^{2^{N}}C_\nu^j$ is bounded by $2^{N}2^{N+1}<4^{N+1}
$, as claimed.

To finish the induction step we repeat this argument to construct
$b_{-j}$.
For $\nu=1,\dots, 2^{N}$ we denote by $C_\nu^{-j}$
the subset of all integers $c$ in
$ [-j4^{N+1},
(1-j) 4^{N+1})$ with the property that
$e_\nu+c\in\cup_{i=1-j}^{j}(b_i+E).$
Again we verify (by repeating the argument above) that
the cardinality of
$\cup_{\nu=1}^{2^{N}}C_\nu^{-j}$ is $ < 4^{N+1}$ and then we may choose
$b_{-j}\in
[-j4^{N+1}, (1-j) 4^{N+1})$ so that $b_{-j}$ does not belong to
$\cup_{\nu=1}^{2^{N}}C_\nu^{-j}$. Then by construction
 the sets $b_{-j}+E, \dots, b_j+E$ are disjoint.
\qed

\bigskip

\noi{\bf Weak type (1,1) and Hardy space estimates.} 
For a countable   set of  multipliers $\{m_{\nu}\}_{\nu\in\cI}$
consider the maximal function given by
$$\fM f(x)= \sup_{\nu \in \cI} |\cF^{-1}[m_\nu \widehat f\, ](x)|.$$
We shall apply   the following result
on maximal functions
which is based on Calder\'on-Zygmund
theory  and essentially proved in \cite{DT}. In what follows $H^1$ denotes the standard Hardy space.

\medskip

\noi {\bf Proposition 3.2.}  \emph{
Suppose that for some positive $\epsilon \le 1$}
$$\sup_{\nu\in \cI }\sup_{k\in \bbZ}
 %% \|\phi m_\nu(2^k \cdot)\|_{L^1_{d+\epsilon}}\le A_0.
%\int
\sup_{x\in \bbR^d} (1+|x|)^{d+\epsilon}
\big|\cF^{-1}[\phi m_\nu(2^k \cdot)](x)\big| \le  A_0
%%\|_{L^1_{d+\epsilon}}\le A_0.
$$
\emph{$\!\!\!$
and suppose that $\fM$ is bounded on $L^q$ (for some $q>1$)
with operator norm $A_1$. Then  $\fM$ is bounded from $H^1$ to $L^1$ with 
operator norm at most  $  C_d (A_0\eps^{-1}+A_1)$; moreover $\fM$ is  of weak type (1,1) 
with the estimate}
$$
\sup_{\alpha>0}\,\alpha\,\meas\big(\{x:|\fM f(x)|>\alpha\}\big)
\le C_{d} (A_0\eps^{-1}+A_1)\|f\|_1.
$$

\noindent{\it Proof.} We prove the weak-type $(1,1) $ bound.
Fix $\alpha>0$.
We use the standard Calder\'on-Zygmund decomposition  (see \cite{Stein})
at level $\beta= (2^{d+1}A_1)^{-1}\alpha$.
Thus we decompose $f=g_\beta+b_\beta$ where $|g_\beta|\le 2^d \beta$ and
 $b_\beta=\sum b_{\beta,Q}$, where
$b_{\beta,Q}$ is supported on $Q$ and has mean value $0$.
 Moreover, if $Q^*$ denotes the $2\sqrt d$-dilate of $Q$  with 
same center,  then the dilated cubes $Q^*$ have bounded overlap and  $$\sum_Q \meas(Q^*)\le C(d)\beta^{-1} \|f\|_1\le
2^{d+1}C(d) A_1 \alpha^{-1}\|f\|_1.
$$

Let $K_{\nu,j}=\cF^{-1}[ \phi(2^{-j}\cdot) m_\nu]$. We argue similarly as in
 Lemma 1 of
\cite{DT} to verify   the following  vector-valued H\"ormander condition
for maximal operators
(see \cite{Zo}):
\begin{equation}
\label{hoermander-condition}
\int_{|x|\ge 2|y|} \sup_{\nu} \sum_j \Big| K_{\nu,j}(x-y)  -
K_{\nu,j}(x)\big| dx
\le  C_{d} \epsilon^{-1} A_0.
\end{equation}

Let $\widetilde K_{\nu,j}=
\cF^{-1}[ \phi m_\nu (2^j\cdot)]$ so that
$K_{\nu,j}(x)= 2^{jd} \widetilde K_{\nu,j}(2^jx)$.
By assumption we have  the pointwise estimate
$$
|\widetilde K_{\nu,j}(x)|+
|\nabla \widetilde K_{\nu,j}(x)|
\lc A_0 (1+|x|)^{-d-\eps},
$$
uniformly  in $\nu$ and $j$. This quickly yields
$$\int_{|x|\ge 2|y|} \sup_{\nu}  \Big| K_{\nu,j}(x-y)  -
K_{\nu,j}(x)\big| dx
\lc A_0 \min\{(2^j|y|)^{-\eps}, 2^j|y|\};
$$
thus after summing in $j$ we obtain
\eqref
{hoermander-condition}.
This inequality  implies in the usual way
$$
\meas\{x\notin \cup Q^* : \fM b_\beta(x)>\alpha/2\} \lc
A_0\eps^{-1}\alpha^{-1} \|f\|_1.
$$
For the contribution of the ``good'' function $g_\beta$ we obtain
\begin{equation}\label{good}
\begin{aligned}
&\meas\{x: \fM g_\beta(x)>\alpha/2\} \le 2^q \alpha^{-q}  \| \fM
g_\beta\|_q^q
\\
&\le 2^q \alpha^{-q}  A_1^q \|g_\beta\|_q^q
\le 2^q \alpha^{-q}   A_1^q (2^d \beta)^{q-1} \|f\|_1
\\
&\le 2 A_1 \alpha^{-1}\|f\|_1.
\end{aligned}
\end{equation}
A combination of these estimates yields the weak-type $(1,1)$ estimate.

Finally for the $H^1-L^1$ bound we use the atomic decomposition of $H^1$ and it suffices to prove the estimate 
$$\|\fM f_Q\|_1\lc (A_0\eps^{-1}+A_1)$$
for functions $f_Q$ supported on a cube $Q$ satisfying $\int f_Q dx=0$ and 
$\|f_Q\|_\infty\le|Q|^{-1}.$ If $Q^*$ denotes the expanded cube,
 then we get  
$$
\|\fM f_Q\|_{L^1(Q^*)}  \lc |Q| ^{\frac 1{q'}} \| \fM f_Q\|_q\lc
 |Q| ^{\frac 1{q'}} A_1\| f_Q\|_q \lc
 A_1.$$
Using the cancellation of the atom we see that
\eqref
{hoermander-condition} implies 
$$\|\fM f_Q\|_{L^1(\bbR^d\setminus Q^*)}\lc A_0\eps^{-1}$$
and combining the two estimates we get the asserted $H^1\to L^1$ estimate.
\qed

\medskip

Note that the hypothesis in the proposition is implied by
$$
\sup_{\nu\in \cI }\sup_{k\in \bbZ}
 \|\phi m_\nu(2^k \cdot)\|_{L^1_{d+\epsilon}}\le A_0.
$$
The result of Dappa and Trebels \cite{DT} mentioned in the introduction corresponds to the special case where $m_\nu= m(t_\nu\cdot)$ and $\{t_\nu\}$ is 
an enumeration of the positive rational numbers.

\section{Results on $L^p$ boundedness}

In this section  $E$ will be a set of integers satisfying
\begin{align}
%\label{mod5}
%k\in E, \, k'\in E &\, \implies \, |k-k'|>4,
%\\
\label{cardinality}
\card (E) &\le 2^N;
\end{align}
here $N$ is a nonnegative integer.
%Also we shall denote by $\fA(E)$ the union of disjoint annuli
%\begin{equation}
%\label{annuli}
%\fA(E)= \bigcup_{k\in E} \{\xi: 2^k\le |\xi|\le 2^{k+1}\}.
%\end{equation}

Let $\psi\in C^\infty_c$ be supported in 
$\{ \xi:1/4<|\xi|< 4\}$ and set $\Psi=\cF^{-1}[\psi]$;
later we shall work with a specific $\psi$ satisfying \eqref{reproducing}.

Let $\chi\in C^\infty_c(\bbR^d)$ so that $\chi$ is radial
 and supported where
$R_0\le |x|\le R_1$ with $1/2<R_0<R_1<2$; moreover assume that $\chi$ 
is positive for $R_0<|x|<R_1$ and that $\sum_{l=-\infty}^\infty 
 \chi(2^{-l}x)=1$ for $x\neq 0$.
Now set $\chi_0(x)= \big(1-\sum_{l>0} \chi(2^{-l}x)\big)$ so that $\chi_0$ is
supported  where $|x|\le 2$. Let 
$\chi_l(x)=\chi(2^{-l}x)$ for $l>0$; then 
$\sum_{l=0}^\infty \chi_l(x)\equiv 1$.

For a function $g$ define by $\delta_t g$ the $L^1$ dilate; {\it i.e.}
$$\delta_t g(x) = t^{-d}g(t^{-1} x)
% = \cF^{-1}[\widehat g(t\cdot)](x)
.$$
For a sequence $H=\{h_k\}_{k\in \bbZ}$ 
of locally integrable functons we then consider the operator

\begin{equation}\label{TEl}\cT_t^{E,l}[H,f]=
\sum_{k\in E} \delta_{2^kt}\big[\Psi*(\chi_l h_k)\big]  * f
\end{equation} 
and the maximal function
$$\cM^{E,l}[H,f]=\sup_{t>0}\big|\cT_t^{E,l}[H,f]|.$$

In \S5 we shall decompose $\cF^{-1}[m(t\cdot)\widehat f]$  in terms of 
operators of the form \eqref{TEl}.

The following $L^q$ bound is favorable when $q\ge N+l$.

\medskip
\noindent{\bf Proposition 4.1.} \emph{ Assuming
\eqref{cardinality} we have for $q\ge 2$}
%%Suppose that $\card (E) \le 2^N$. Then} 
$$
\|\cM^{E,l} [H,f]\|_q \le C 
  \,q \, 4^{N/q}\,  (1+l)  2^{l/q}\, \|H\|_{\ell^\infty(L^1)}
\|f\|_q.
$$

\medskip

%For the proof of the proposition we shall need the following

\noindent{\it Proof.} 
In what follows we shall use the
notation $A\lc B$ to indicate an inequality
$A\le C B$ where $C$ may only depend on $d$
(and  not on $q$ or other parameters).

Define $$g_{k,l}(\xi)= \psi(\xi) \widehat {h_k\chi_l}(\xi)$$
and 
\begin{equation}\label{defofml}
m^l(\xi)\equiv m^{l,E}(\xi) 
=\sum_{k\in E} g_{k,l}(2^{-k}\xi);
\end{equation}
then $$\cT^{E,l}_t [H,f]= \cF^{-1}[m^l(t\cdot)\widehat f].
%%=\sum_{k\in E} g_{k,l}(2^{-k}t\xi) \widehat f(\xi).
$$
Also note that
%the formula
\begin{equation}\label{ml-deriv}
\partial_s \big[m^l(s\xi)\big]= \sum_{k\in E} \tg_{k,l}(s2^{-k}\xi)
\end{equation}
where
\begin{equation}\label{tgkl}
\tg_{k,l}(\xi) = \psi(\xi)\inn{\xi}{\nabla} \widehat {\chi_l h_k}(\xi)+
\widehat {\chi_l h_k}(\xi)
\inn{\xi}{\nabla}  \psi(\xi).
\end{equation}

Now apply Lemma 3.1 for the set $E$, and let 
$b_j$ be as in Lemma 3.1 (ii).
By (iii) of Lemma 3.1 we may write
$$\sup_{t>0}|\cF^{-1}[m^l(t\cdot) \widehat f \,]|=
\sup_{|n|\le 4^{N+1}}\sup_{j\in \bbZ}\sup_{1\le s\le 2}
\big|\cF^{-1}[m^l(2^{-b_{j}+n}s\cdot)\widehat f \,]\big|.
$$
Now one  replaces the supremum in $n$ and $j$ by $\ell^q$ norms,
takes the $L^q$ norms, then interchanges  the order of summation and
integration. This yields for
$\cM^{E,l}[H,f]\equiv \cM_{m^l} f$
the estimate 
\begin{equation}
\label{lemma32}
\big\|\cM_{m^l}
f\big\|_q\le 4^{\frac{(N+2)}{q}}
\sup_{|n|\le 4^{N+1}}\Big(\sum_j\Big\| \sup_{1\le s\le 2}\big|
\cF^{-1}[m^l(2^{-b_{j}+n}s\cdot)\widehat f \, ]\big|
\Big\|_q^q\Big)^{\frac1q}.
\end{equation}

Thus it remains to verify that for
$|n|\le 4^{N+1}$
\begin{equation}
\label{lqbound}
\Big(\sum_j\Big\| \sup_{1\le s\le 2}\big|
\cF^{-1}[m^l(2^{-b_{j}+n}s\cdot)\widehat f \,]\big|
\Big\|_q^q\Big)^{\frac{1}{q}}
\,\lc \,q(1+l) 2^{\frac lq}\|H\|_{\ell^\infty(L^1)}
\|f\|_q.
\end{equation}
In what follows we
may assume  that $n=0$ since the general case follows by scaling.

To estimate the supremum in $s$ it is standard to
use the elementary inequality
$$
|F(s)|^q\le |F(1)|^q + q
\Big(\int_1^s |F(\sigma)|^q d\sigma\Big)^{\frac{q-1}{q }}
\Big(\int_1^s |F'(\sigma)|^q d\sigma\Big)^{\frac{1}{q }}
$$
which is obtained by applying the
fundamental theorem of calculus to
$|F|^q$ and   H\"older's inequality.

Taking $L^q$ norms and applying  H\"older's
inequality twice  yields
\begin{equation}
\label{1-2-3}
\begin{aligned}
\sum_j\big\| \sup_{1\le s\le 2}\big|
\cF^{-1}[&m^l(2^{-b_{j}}s\cdot)\widehat f\, ]\big|
\big\|_q^q  \lc
\sum_{j\in \bbZ}  \big\|III^l_{j}\big\|_q^q
\,
\\
+\, q  &\,
\Big(\int_1^2 \sum_{j\in \bbZ} \big\|I^l_{j}(s)
\big\|_q^q ds\Big)^{\frac{q-1}{q }}
\Big(\int_1^2 \sum_{j\in \bbZ}
\big\|II^l_{j}(s)\big\|_q^q
ds\Big)^{\frac{1}{q }}\, ,
\end{aligned}
\end{equation}
where
\begin{align*}
I_{j}^l(s)
&=\cF^{-1}\big[
m^l(2^{-b_{j}}s\cdot)
\widehat f\, \big]
\\
II_{j}^l(s)&=\cF^{-1}\big[
\partial_s\big(m^{l}(2^{-b_{j}}\cdot)\big) \widehat f\, \big]
\\
III_{j}^l&=\cF^{-1}\big[
m^{l} (2^{-b_{j}}\cdot) \widehat f\, \big].
\end{align*}

Next,  we interchange the  $j$-summations and integrations in \eqref{1-2-3}
 and use the imbedding
of $\ell^2$ into $\ell^q$. This
yields

\begin{equation}
\label{threeterms}
\begin{aligned}
&\Big(\sum_j\big\| \sup_{1\le s\le 2}|
\cF^{-1}[m^{l}(2^{-b_{j}}s\cdot)\widehat f \, ]|
\big\|_q^q \Big)^{1/q}
\\&\,\,\lc
\Big(\int_1^2
\Big\|\Big(\sum_{j\in \bbZ} |I^l_{j}(s)|^2\Big)^{1/2}
\Big\|_q^q ds\Big)^{\frac{q-1}{q^2 }}
%%\\&\qquad\times
\Big(\int_1^2\Big\|\Big(\sum_{j\in \bbZ}
|II^l_{j}(s)|^2\Big)^{1/2}\Big\|_q^q
ds\Big)^{\frac{1}{q^2}}
\\&\qquad +
\Big\|\Big(
\sum_{j\in \bbZ}  |III^l_{j}|^2\Big)^{1/2}\Big\|_q.
\end{aligned}
%\end{multline*}
\end{equation}

In order to estimate these terms we need
 the following estimates for vector valued singular integrals.

\medskip

\noindent{{\bf Sublemma.}} \emph{ For $2\le q<   \infty$ we have  }
%$|n|\le  4^{N+1}$,
\begin{equation}
\label{sublemma1}
\Big\|\Big(\sum_{j\in \bbZ}\Big|\sum_{k\in E} \cF^{-1}[
g_{k,l}(2^{-b_{j}-k}\cdot) \widehat f\, ] \Big|^2\Big)^{1/2}\Big\|_q
\lc q (1+l) \|H\|_{\ell^\infty(L^1)}
\|f\|_q
\end{equation}
\emph{and}
\begin{equation}
\label{sublemma2}
\Big\|\Big(\sum_{j\in \bbZ}\Big|\sum_{k\in E} \cF^{-1}[
\tg_{k,l}(2^{-b_{j}-k}\cdot) \widehat f \, ] \Big|^2\Big)^{1/2}\Big\|_q
\lc q (1+l) 2^l \|H\|_{\ell^\infty(L^1)}
\|f\|_q.
\end{equation}

\medskip

\noindent {\it Sketch of  Proof.}
By duality  \eqref{sublemma1} for $2\le q<\infty$  is equivalent to
\begin{multline}\label{sublemma-duality}
\Big\|\Big(\sum_{j\in \bbZ}\Big|\sum_{k\in E} \cF^{-1}[
g_{k,l}(2^{-b_{j}-k}\cdot) \widehat {f_j}] \Big|^2\Big)^{1/2}\Big\|_p
\\
\lc p'
(1+l)^{-1+2/p} \|H\|_{\ell^\infty(L^1)} 
\Big\|\Big(\sum_j|f_j|^2\Big)^{1/2}\Big\|_p
\end{multline}
for $1<p\le 2$, $p'=p/(p-1)$. For $p=2$ this (and in fact a slightly better) 
bound  follows from the 
essential disjointness of the supports of $g_{k,l}$ and the estimate
$$\|g_{k,l}\|_\infty \le C \|\widehat {h_k}\|_\infty\le C \|H\|_{\ell^\infty(L^1)}.$$
For $1<p\le 2$ the inequality \eqref{sublemma-duality}
 follows from the weak type bound
\begin{multline}\label{weaktype}
\meas\Big(\Big\{x:
\Big(\sum_{j\in \bbZ}\Big|\sum_{k\in E} \cF^{-1}[
g_{k,l}(2^{-b_{j}-k}\cdot) \widehat {f_j}] \Big|^2\Big)^{1/2}>\la
\Big\}\Big)
\\
\lc C \la^{-1} 
(1+l) \|H\|_{\ell^\infty(L^1)} 
\Big\|\Big(\sum_j|f_j|^2\Big)^{1/2}\Big\|_1
\end{multline} 
where in the interpolation we have to take into account the 
behavior of the constants in the Marcinkiewicz interpolation theorem
(see e.g.  \cite{G}, p. 33). 

The  weak type estimate follows by standard arguments 
in Calder\'on-Zygmund theory from
the inequality 
\begin{multline*}
\int\limits_{|x|>2|y|}
\Big(\sum_{j}\Big|\sum_{k\in E}\Big(
\cF^{-1}[g_{k,l}(2^{-b_j-k}\cdot)](x-y)\\
-
\cF^{-1}[g_{k,l}(2^{-b_j-k}\cdot)](x)\Big)\Big|^2\Big)^{1/2} dx 
\lc (1+l)
\|H\|_{\ell^\infty(L_1)}
\end{multline*}
which, since the sets $\{b_j+E\}_{j\in \bbZ}$ are disjoint, is 
quickly derived from the inequalities
\begin{multline}
\int_{|x|>2|y|}\label{estimateforpieces}
\big|
\cF^{-1}[g_{k,l}(2^{-M}\cdot)](x-y)
-
\cF^{-1}[g_{k,l}(2^{-M}\cdot)](x)\big| dx
\\
\lc \|h_k\|_1 \times 
\begin{cases}
2^M |y| \quad&\text{ if } 2^M|y|\le 1
\\
1 \quad&\text{ if } 1\le  2^M|y|\le 2^{2l}
\\
(2^M|y|)^{-1} \quad&\text{ if } 2^M|y|\ge 2^{2l}
\end{cases}.
\end{multline}
Summing in $M$ yields a blowup of order  $O(1+l)$. The bound 
\eqref{estimateforpieces} is straightforward given the localization of 
$\chi_l$ and the decay of the Schwartz-function  $\Psi$.
This finishes the proof of
\eqref{sublemma1}.

In order to verify \eqref{sublemma2} we note from \eqref{tgkl} that 
\begin{multline*}
\cF^{-1}[\tg_{k,l}](x)=
c_1\sum_{i=1}^d 
\int \partial_{x_i} \Psi(x-y) y_i\chi_l(y) h_k(y)
dy\\+
c_2\sum_{i=1}^d 
\int x_i\partial_{x_i}\Psi (x-y)\chi_l(y) h_k(y) dy.
\end{multline*}
Here the second  term has the same quantitative properties as 
$\cF^{-1}[g_{k,l}]$ while the first has similar estimates as $ 2^l \cF^{-1}[g_{k,l}]$. Thus the above arguments show \eqref{sublemma2} as well.
\qed

\medskip

\noi{\it Proof of Proposition 4.1, cont.}
For fixed $s$ we may perform the scaling $\xi\to  s^{-1}\xi$ in
\eqref{threeterms}; this puts us in the position to apply the sublemma. We
then see that the right hand side of \eqref{threeterms} can be estimated by a constant times
$$ q(1+l) \|H\|_{\ell^\infty(L^1)} (1+2^{l/q})\|f\|_q
$$
which implies the desired  bound \eqref{lqbound}. \qed

\medskip

\noindent{\bf Corollary 4.2.} \emph{
Assuming \eqref{cardinality}  we have}
% Suppose that
%$\card (E) \le 2^N$. Then for $1<p<\infty$
$$
\|\cM^{E,l} [H,f]\|_{BMO}
\le C
(N+l)(1+l) \|H\|_{\ell^\infty(L^1)}  \|f\|_{L^\infty}.
$$

\medskip

\noi{\it Proof.}
Let $Q$ be a cube in $\bbR^d$ with center $x_Q$ and let $Q^*$ be
the  $2\sqrt d$-dilate with same center.
By the definition of $BMO$ and H\"older's inequality 
the assertion follows from 
\begin{equation}\label{bmo1}
\Big(\frac 1{|Q|}\int_Q |\cM^{E,l}[ H, f\chi_{Q^*}]|^{N+l} dx \Big)^{
\frac {1}{N+l}}
\lc(N+l)(1+l) \|H\|_{\ell^\infty(L^1)}  \|f\|_{L^\infty}
\end{equation}
and, with $m^l$ as in \eqref{defofml},
\begin{multline}\label{bmo2}
\sup_{x\in Q} \sup_{t>0}\int_{\bbR^d\setminus Q^*} \big|\cF^{-1}[m^l(t\cdot)](x-y)
-\cF^{-1}[m^l(t\cdot)](x_Q-y)\big| dy \\
\lc
(1+l) \|H\|_{\ell^\infty(L^1)}  \|f\|_{L^\infty}.
\end{multline}
The left hand side  of \eqref{bmo1} 
is bounded by the $L^{N+l}$ operator norm of $\cM^{E,l}(H, \cdot)$
times $|Q|^{-1/(N+l)}\|f\chi_{Q^*}\|_{N+l}$ and 
\eqref{bmo1} follows from Proposition 4.1. Inequality \eqref{bmo2} is deduced from the estimates in the Sublemma.\qed

\medskip

\noindent{\bf Proposition 4.3.} 
\emph{Suppose that \eqref{cardinality} holds. Then
 for $1<p<\infty$, $1/p+1/p'=1$,} 
$$
\|\cM^{E,l} [H,f]\|_p
\le C_p (N+l)(1+l)
 2^{ld/p} \|H\|_{\ell^\infty(L^{p'})}  \|f\|_p.
$$
\emph{
Moreover  the operator $f\mapsto\cM^{E,l} [H,f]$
is bounded from $H^1$ to $L^1$, and of weak type $(1,1)$ 
with operator bound $C N(1+l)
 2^{ld} \|H\|_{\ell^\infty(L^{\infty})}$.}

\medskip

\noi{\it Proof. } 
Let $\tilde \chi_l$ be the characteristic function of the annulus
$\{x:2^{l-3}\le |x|\le 2^{l+3}\}$, for $l\ge 1$, and let $\tilde 
\chi_0$ be the characteristic function of the ball 
$\{x:|x|\le 8\}$. Let 
$h_k^l= h_k \tilde \chi_l$ and $H^l=\{h_k^l\}_{k\in \bbZ}$. 
Then observe that by the support property of $\chi_l$ 
(in the definition of $\cT^{E,l}$) we have
$\cT^{E,l}[H,f]=\cT^{E,l}[H^l, f]$.
Proposition 4.1 yields the estimate
$$
\|\cM^{E,l} [H,f]\|_q = \|\cM^{E,l} [H^l,f]\|_q
 \le C q(1+l) \|H^l\|_{\ell^\infty(L^1)}, \quad q\ge N+l
$$
This implies the assertion 
%(or a slightly better version) 
 for $ N+l\le p<\infty$.

To prove the  Hardy space estimate 
we apply Proposition 3.2 in conjunction with 
Proposition 4.1 (for $H^l$ and $q=N+l$) and we obtain
\begin{multline*}\|\cM^{E,l}[H,f]\|_1\lc 
\\
\big((N+l)(1+l)\|H^l\|_{\ell^\infty(L^1)}
+(1+l)\eps^{-1} 2^{l\epsilon} 2^{ld}\|H^l\|_{\ell^\infty(L^\infty)}\big) 
\|f\|_{H^1}.
\end{multline*}
The asserted bound follows if we choose $\eps=(1+l)^{-1}$.
 The  weak type $(1,1)$ bound follows similarly.

We may now use the complex method for bilinear operators 
(which is a variant of Stein's theorem for analytic families, see \cite{BL}, 
\S 4.4) together with the interpolation formula $[H_1, L^{p_1}]_\vartheta=L^p$,
 for 
$(1-\vth)+\vth/p_1=1/p$, see \cite {FS}.
Now define $s=p'/p_1'$ so that $\vartheta=(1-\vartheta)/\infty+\vartheta/1=1/s$.
We then obtain the estimate
$$
\|\cM^{E,l} [H,f]\|_p \le C_p (1+l)(N+l)2^{ld/s'}
\|H\|_{\ell^\infty(L^s)}\|f\|_p
$$
but since 
$\cM^{E,l} [H,f] = \cM^{E,l} [H^l,f]$ we may replace $H$ with $H^l$ on the right hand side of this inequality. 
Note that $s<p'$ and thus
$$2^{ld/s'} \|H^l\|_{\ell^\infty(L^s)}
\lc 2^{ld/p} \|H^l\|_{\ell^\infty(L^{p'})},$$ 
by H\"older's
inequality. This yields the asserted bound.\qed
% Marcinkiewicz interpolation theorem (see e.g.  \cite{G}, p. 33). 
%This yields
%the $L^p$ inequality with constant
%$\lc (1+l) (p'+ p/(2(N+l)-p))^{1/p} 
%\|H^l\|_{L^\infty(\ell^\infty)} $ which 
%implies the assertion for $1<p\le N+l$. \qed

\medskip

\noi{\it Remark}. One could also analytically 
interpolate the $H^1\to L^1$ estimate with the $L^\infty\to BMO$ estimate 
of Corollary 4.2 using the formula $[H_1,BMO]_\theta= L^p$, $\theta=1/p'$.
This formula follows  from the result in \cite{FS} for 
$[H_1, L^{p_1}]_{\vartheta_1}$ and  
$[L^{p_0}, BMO]_{\vartheta_2}$, $1<p_0<p_1<\infty$,   by 
Wolff's four space reiteration theorem, see \cite{W}.

\medskip

\section{Conclusion: Proof of Theorem 1.2}
We only have to prove the $L^p$ estimates for $p>1$ since the asserted weak
type $(1,1)$ bound is then a consequence of Proposition 3.2.

We need to decompose $m$ in terms of
the rearrangement function $\omega^*$.
Let
$$
%\begin{equation}
E_0= \{k\in \bbZ: \omega^*(2)<|\omega(k)|\le \omega^*(0)\}
$$
%\end{equation}
and for
$j=1,2,\dots$  let
%\begin{equation}
$$E_j= \{k\in \bbZ: \omega^*(2^{2^j})<|\omega(k)|\le \omega^*(2^{2^{j-1}})\}.$$ 
%\end{equation}
%We further split $E_N$ into five disjoint subsets,
%$E_N=\cup_{i=1}^5 E_{N,i}$ where each $E_{N,i}$ satisfies \eqref{mod5}.
%Let
%\begin{equation}
%\label{mN}

As in the introduction  let $\phi\in C^\infty_c$ be   supported in
$\{ \xi:1/2<|\xi|< 2\}$ so that $\phi(\xi)\neq 0$ for
$2^{-1/2}\le |\xi|\le 2^{1/2}$. Set $\psi(\xi)=\overline{\phi(\xi)}
\big(\sum_{j\in \bbZ}|\phi(2^{-j}\xi)|^2\big)^{-1}$, then $\psi$ is smooth and
we have
\begin{align}
\label{reproducing}
\sum_{k\in \bbZ} \psi(2^{-k}\xi) \phi(2^{-k}\xi)  =1,\quad \xi\neq 0.
\end{align}

Let
$$m_{j}(\xi)=\sum_{k\in E_{j}}
   \psi(2^{-k}\xi)\phi(2^{-k}\xi)
m(\xi)$$
%\end{equation}
then
$m= \sum_{j=0}^\infty m_j$; here we use that
$\omega^*(k)\to 0$ as $k\to \infty$.

%Now let $1<p<\infty$ and assume that \eqref{thm12assumption}
% holds.

%We use for a Schwartz function $\zeta$ the elementary inequality
%$$
% \Big(\int \big|G*\zeta(x)|^{p'}(1+|x|)^{\alpha p'} 
%dx\Big)^{1/p'} \lc
% \Big(\int \big|G(x)|^{p'}(1+|x|)^{\alpha p'} 
%dx\Big)^{1/p'}.
%$$
%If we apply this to 
%$\zeta=\sum_{i=-2}^2\cF^{-1}[\phi(2^{-i}\cdot) \psi(2^{-i}\cdot)]$
%then  we see that also
%$$
% \Big(\int \big|\cF^{-1}[\phi m_j(2^k\cdot)]\big|^{p'}(1+|x|)^{\alpha p'} 
%dx\Big)^{1/p'} \le \omega^*(2^{2^{j-1}})
%$$

%We have  for $k\in E_\ell$
%\begin{align*}
%\|\phi m_\ell(2^k\cdot)\|_{L^r_\alpha}&\le
%\sum_{j=k-2}^{k+2} \|\phi m(2^k\cdot)\phi(2^{k-j}\cdot)\psi(2^{k-j}\cdot)
%\|_{L^r_\alpha}
%\\& \lc\|\phi m(2^k\cdot)\|_{L^r_\alpha}\lc
%\omega^*(2^{2^{\ell-1}}).
%\end{align*}
%Note from the definition of  the rearrangement function that
%\begin{equation}\label{cardN}
%$\card(E_\ell)\le 2^{2^\ell}$.
% \end{equation}

We now decompose $\cF^{-1}[\phi m(2^k\cdot)]$ using the dyadic cutoff functions
$\chi_l$.
Define
$h_{k}^{j,0}(x)=\cF^{-1}[\phi m(2^k\cdot)](x)$ if $|x|\le 4$ and $k\in E_j$ 
and $h_{k}^{j,0}(x)=0$, if $|x|>4$ or $k\notin E_j$. Moreover 
for $l>0$ let $h_{k}^{j,l}(x)= 
\cF^{-1}[\phi m(2^k\cdot)](x)$ if $2^{l-4}\le |x|\le 2^{l+ 4}$ and $k\in E_j$
 and $h_{k}^{j,l}(x)=0$, if $|x|\notin [2^{l-4}, 2^{l+4}]$ or $k\notin E_j$. 
Then
$$
\cF^{-1}[\phi m(2^k\cdot)]= \sum_{j=0}^\infty \sum_{l=0}^\infty h_{k}^{j,l}\chi_l
$$
and with  $H^{j,l}=\{h_k^{j,l}\}_{k\in \bbZ}$ , and $\Psi=\cF^{-1}\psi$, 
 we may write 
\begin{align*}
\cF^{-1}[m_j(t\cdot)\widehat f]
&=\sum_{k\in E_j} \delta_{2^k t}\big(\Psi*\sum_{l=0}^\infty(h_{k}^{j,l}\chi_l)
\big)
\\
&= \sum_{l=0}^\infty \cT_t^{E_j,l}[H^{j,l}, f].
\end{align*}

The assumption \eqref{thm12assumption} implies that 
$$\|H^{j,l}\|_{\ell^\infty(L^{p'})}\lc 2^{-l\alpha} \omega^*(2^{2^{j-1}})$$
for $j\ge 1$ (and a similar estimate with $\omega^*(0)$ for $j=0$).

Note from the definition of the rearrangement function that 
 $$\card(E_j)\le 2^{2^j}.$$  Thus we obtain by Proposition 4.3
$$
\big\|\cM^{E_j,l}[H^{j,l}, f]\big\|_p \lc (1+l)^2 2^{-l(\alpha-d/p')}
2^j\omega^*( 2^{2^{j-1}}), \quad j\ge 1,
$$
and a similar estimate with $2^j\omega^*( 2^{2^{j-1}})$ replaced by 
$\omega^*(0)$ when $j=0$.
By our assumption $\alpha>d/p'$ we therefore  get
\begin{align*}
\|\cM_{m} f\|_p
&\le \sum_{j=0}^\infty \|\cM_{m_j} f\|_p\le \sum_{l=0}^\infty\sum_{j=0}^\infty \|\cM^{E_j,l}[H^{j,l},f]\|_p
\\ &\le C_p \sum_{l=0}^\infty (1+l)^2 2^{-l(\alpha-d/p')}
     \Big[\omega^*(0)+\sum_{j=1}^\infty 2^j
\omega^*(2^{2^{j-1}})\Big]\|f\|_p
\\
&\lc C_p'\Big[\omega^*(0)+\sum_{n=2}^\infty \frac{\omega^*(n)}{n}\Big]\,\|f\|_p.
\end{align*}
This concludes the proof of Theorem 1.2.\qed

\medskip

\section{A sketch of the proof of Theorem 1.4.}
The proof  is  simpler than  the proof of Theorem 1.2 but relies
on the same idea. Write $T_\nu f=\cF^{-1}[m_\nu \widehat f]$.
Assume that
\eqref{thm14assumption}
holds with $p'=1$, $\alpha>0$. Then 
for $10\alpha^{-1}<q<\infty$ the operator $T_\nu$ is bounded on $L^q$ with
 operator norm  $O(q)$,
uniformly in $\nu$. We replace $\ell^\infty$ norms by $\ell^q$
norms and estimate for those $q$ 
\begin{align*}
\|\fM_n f\|_q& \le \Big\|\big(\sum_{\nu=1}^n |T_\nu
f|^q\big)^{1/q}\Big\|_q = \Big(\sum_{\nu=1}^n \|T_\nu
f\|_q^q\Big)^{1/q}\le C_\alpha  n^{1/q} q\|f\|_p.
\end{align*}
This yields the desired result for $q\ge \log(n+1) \ge 10\alpha^{-1}$
 since $n^{1/\log(n+1)}$
is bounded as $n\to \infty$. Arguing as in the proof of Corollary 4.2 one also gets an $L^\infty\to BMO$ estimate with bound $O(\log (n+1))$.
Finally, under the analogue of
\eqref{thm14assumption} for $p=1$ we derive an $H^1\to L^1$ estimate (using Proposition 3.2) and by an interpolation argument as used  in Proposition 4.3
we may derive the asserted  result for $L^p$ boundedness.\qed

\section{Remarks and open problems}

7.1. The main  open problem is to  completely close
 the gap in terms of the power of logarithms in \eqref{example} and
\eqref{logpowerone}. In particular it should be interesting to
know assuming \eqref{sob-Lr} for which $s>1$ the condition
$$\omega^*(0)+\Big(\sum_{l= 1}^\infty \frac{[\omega^*(l)]^s}{l}\Big)^{1/s}<\infty$$
implies $L^p$ boundedness of $\cM_m$.  We have shown that $s=1$ is sufficient,
but $s>2$ is not.

Similarly, in \eqref{lognbound} it would be interesting to
investigate
 whether the bound
$O(\log n)$ can be replaced by $O(\log^{1/s}\! n)$ for suitable
$s\le 2$.

\medskip

7.2. If $d=1$ the  assumptions
\eqref{thm14assumption} on the kernels (which are differentiability 
assumptions on the  multiplier) are essentially sharp; this is seen by 
examining
the multipliers $e^{i\xi}|\xi|^{-\alpha}\sum_{k=1}^N \phi(2^{-k}\xi)$ 
for suitable $\phi$.

\medskip

7.3. If $m$ is radial, $m(\xi)=g(|\xi|)$,  then the space
$L^1_{d+\eps}$ may be replaced by
$L^1_{(d+1)/2+\eps}$ in the weak type (1,1) estimate (see \cite{DT}),
of Corollary 1.3.
By analytic interpolation  one  obtain $L^p$ boundedness for $1\le p\le 2$
under the conditions \eqref{sob-Lr}, \eqref{omega}  with
 with $r=\min\{p,2\}$,  $\alpha>d/2+(1/r-1/2)$.

\medskip

7.4.
 If  \eqref{omega} is replaced by a stronger decay assumption
then much weaker smoothness assumptions suffice, as demonstrated in
\cite{C}, \cite{DT} under the assumption $\omega\in \ell^2$.
 Various intermediate
estimates can be derived by analytic interpolation.  It should be
interesting to  obtain in higher dimensions
 the minimal smoothness assumption requiring only
the decay in \eqref{omega}. The same question can be formulated for the 
dyadic maximal operators.

\end{document}